\documentclass{amsart}
\usepackage{graphicx}

\newtheorem{defn}{Definition}
\newtheorem{rem}{Remark}
\newtheorem{exa}{Example}

\usepackage{color}
\definecolor{brown}{rgb}{0.8,0.6,0.3}
\definecolor{dgreen}{rgb}{0.2,0.4,0.3}

\usepackage[pdftex,pdfauthor={Richard J. Mathar}, colorlinks=true,urlcolor=blue,linkcolor=red,citecolor=brown,
pdfkeywords={{Triple binomial product}, {A361608}, {OEIS}},pdftitle={Sum of Triple Product of Binomial Coefficients},bookmarksopen=true,pdfpagelayout=OneColumn]{hyperref}

\begin{document}

\title[Sum of Triple Product of Binomial Coefficients]{On an Alternating Double Sum of a Triple Product of Aerated Binomial Coefficients}

\author{Richard J. Mathar}
\urladdr{https://www.mpia-hd.mpg.de/~mathar}
\address{Max-Planck Institute of Astronomy, K\"onigstuhl 17, 69117 Heidelberg, Germany}

\date{\today}

\begin{abstract}
The double sum $\sum_{j=0}^m
\sum_{i=0}^j
(-)^{j-i}
\binom{m}{j} \binom{j}{i}\binom{j+k+qi}{j+k}$
with free nonnegative integer parameters $k$ and $q$
is rewritten as hypergeometric series.
Efficient formulas to generate the C-finite ordinary generating
functions are presented.
\end{abstract}

\maketitle

\section{Hypergeometric reduction}
The theme of this manuscript is 
the double sum
\begin{defn}
(main sequence)
\begin{equation}
a_{k,q}(m) \equiv\sum_{j=0}^m 
\sum_{i=0}^j 
(-)^{j-i}
\binom{m}{j} \binom{j}{i}\binom{j+k+qi}{j+k} 
\label{eq.defa}
\end{equation}
\end{defn}
for fixed nonnegative integers $k,q$ and nonnegative $m$.
This investigation was triggered by a problem in Project Euler \cite[Prob. 831]{ProjEuler}.
It is a binomial transform of a sequence defined as
\begin{defn}
(intermediate sequence)
\begin{equation}
b_{k,q}(j) \equiv
\sum_{i=0}^j 
(-)^i
\binom{j}{i}\binom{j+k+qi}{j+k} 
.
\label{eq.defb}
\end{equation}
\end{defn}

\begin{exa}
The case $k=0$, $q=1$ is easily reduced with \cite[(3.9)]{RoyAMM94}
\begin{equation}
a_{0,1}(m) \equiv
\sum_{j=0}^m 
(-)^j
\binom{m}{j} 
\sum_{i=0}^j 
(-)^i
\binom{j}{i}
\binom{j+i}{j}
=
\sum_{j=0}^m 
(-)^j
\binom{m}{j} 
(-)^j\binom{j}{0}
=
\sum_{j=0}^m 
\binom{m}{j}  =2^m.
\end{equation} 
\end{exa}

\begin{exa}
There is only one binomial product term if $m=0$:
\begin{equation}
a_{k,q}(0)=1.
\end{equation}
\end{exa}
\begin{exa}
There are three binomial product terms if $m=1$, two cancel:
\begin{equation}
a_{k,q}(1)
=
\binom{1+k+q}{1+k}
.
\end{equation}
\end{exa}

\section{Resumation along Diagonal/Difference of Index Pairs}

\subsection{Reduction to a Sum over Products of Two Binomials}
The region of the inner sum in \eqref{eq.defa} is $0\le i\le j$ which
may as well be rewritten as $i\le j\le m$ \cite{ChoiCKM18},
\begin{multline}
a_{k,q}(m)\equiv \sum_{i=0}^m \sum_{j=i}^m (-)^{j-i}\binom{m}{j} \binom{j}{i}\binom{j+k+qi}{j+k} 
\\
=\sum_{i=0}^m \sum_{j=i}^m (-)^{j-i}\frac{\Gamma(m+1)\Gamma(j+k+qi+1)}{\Gamma(m-j+1)\Gamma(i+1)\Gamma(j-i+1)\Gamma(1+j+k)\Gamma(1+qi)}.
\end{multline}
If the summation is changed to run over the difference $l=j-i$ 
first, this sum can be written as a terminating Gaussian Hypergeometric Function:
\begin{multline}
a_{k,q}(m)=\sum_{i=0}^m \sum_{l=0}^{m-i} (-)^l\frac{\Gamma(m+1)\Gamma(l+k+(q+1)i+1)}{\Gamma(m-l-i+1)\Gamma(i+1)\Gamma(l+1)\Gamma(1+l+i+k)\Gamma(1+qi)}
\\
=\sum_{i=0}^m 
\frac{\Gamma(m+1)\Gamma(k+(q+1)i+1)}{\Gamma(i+1)\Gamma(1+qi)\Gamma(1+i+k)}
\sum_{l=0}^{m-i} (-)^l\frac{(k+(q+1)i+1)_l}{\Gamma(m-l-i+1)l!(1+i+k)_l}
\\
=\sum_{i=0}^m 
\frac{\Gamma(m+1)\Gamma(k+(q+1)i+1)}{\Gamma(i+1)\Gamma(1+qi)\Gamma(1+i+k)}
\sum_{l=0}^{m-i} \frac{(k+(q+1)i+1)_l(i-m)_l}{\Gamma(m-i+1)l!(1+i+k)_l}
\\
=\sum_{i=0}^m 
\frac{\Gamma(m+1)\Gamma(k+(q+1)i+1)}{\Gamma(i+1)\Gamma(1+qi)\Gamma(1+i+k)\Gamma(m-i+1)}
{}_2F_1\left(\begin{array}{c} k+(q+1)i+1, i-m \\ 1+i+k \end{array}\mid  1\right)
.
\end{multline}
This $_2F_1$ can be written as a $\Gamma$-ratio
\cite[15.4.24]{DLMF}\cite[15.1.20]{AS}
\begin{multline}
a_{k,q}(m)=\sum_{i=0}^m 
\binom{m}{i}
\frac{\Gamma(k+(q+1)i+1)}{\Gamma(1+qi)\Gamma(1+i+k)}
\frac{(1+i+k-(k+(q+1)i+1))_{m-i}}{(1+i+k)_{m-i}}
\\
=\sum_{i=0}^m 
\binom{m}{i}
\frac{\Gamma(k+(q+1)i+1)}{\Gamma(1+qi)\Gamma(1+i+k)}
\frac{(-iq)_{m-i}}{(1+i+k)_{m-i}}.
\end{multline}
The Pochhammer symbol of the negative argument is lifted
with \cite{SlaterHyp}
\begin{equation}
(-n)_k = (-1)^k\Gamma(n+1)/\Gamma(n+1-k)
\label{eq.negarg}
\end{equation}
and the double sum \eqref{eq.defa} 
becomes
a single sum:
\begin{multline}
a_{k,q}(m)=\sum_{i=0}^m 
\binom{m}{i}
\frac{\Gamma(k+(q+1)i+1)}{\Gamma(1+qi)\Gamma(1+i+k)}
\frac{(-)^{m-i}\Gamma(iq+1)\Gamma(1+i+k)}{\Gamma(iq+1-(m-i))\Gamma(1+i+k+m-i)}
\\
=\sum_{i=0}^m 
\binom{m}{i}
(-)^{m-i}
\frac{\Gamma(k+(q+1)i+1)}{\Gamma(1+i+k)}
\frac{\Gamma(1+i+k)}{\Gamma(iq+1-(m-i))\Gamma(1+i+k+m-i)}
\\
=\sum_{i=0}^m 
(-)^{m-i}
\binom{m}{i}
\Gamma(k+(q+1)i+1)
\frac{1}{\Gamma((q+1)i-m)\Gamma(1+k+m)}
\\
=\sum_{i=0}^m 
(-)^{m-i}
\binom{m}{i}
\binom{k+(q+1)i}{k+m}
.
\label{eq.2sum}
\end{multline}

\subsection{Transformation to Terminating Generalized Hypergeometric Series}
The previous equation equals a terminating Generalized Hypergeometric Series:
\begin{multline}
a_{k,q}(m) =\sum_{i=0}^m 
(-)^{m-i}
\frac{\Gamma(m+1)\Gamma(k+(q+1)i+1)}{\Gamma(i+1)\Gamma(m-i+1)}
\frac{1}{\Gamma((q+1)i+1-m)\Gamma(1+k+m)}
\\
=
\frac{\Gamma(m+1)}{\Gamma(1+k+m)}
\sum_{i=0}^m 
(-)^{m-i}
\frac{\Gamma(1+k+(q+1)i)}{\Gamma(i+1)\Gamma(m-i+1)}
\frac{1}{\Gamma(1-m+(q+1)i)}
\\
=
\frac{\Gamma(m+1)}{\Gamma(1+k+m)}
\sum_{i=0}^m 
(-)^{m-i}
\frac{\Gamma(1+k) (1+k)_{i(q+1)}}{\Gamma(i+1)\Gamma(m-i+1)}
\frac{1}{\Gamma(1-m) (1-m)_{i(q+1)}}
.
\end{multline}
Where poles appear the denominator, the contribution of the
$i$-term to the sum is $1/\infty=0$.
The Pochhammer analog to the product formula
of the $\Gamma$-function is \cite[(2.4.5.2.)]{SlaterHyp}\cite{DriverETNA25,RoyAMM94}
\begin{equation}
(a)_{qr}=(a/q)_r (\frac{a+1}{q})_r\cdots (\frac{a+q-1}{q})_r  q^{qr}
\label{eq.pmulti},
\end{equation}
so
\begin{equation}
a_{k,q}(m)
=
\frac{\Gamma(m+1)}{\Gamma(1+k+m)}
\sum_{i=0}^m 
(-)^{m-i}
\frac{\Gamma(1+k) \prod_{l=0}^q (\frac{k+1+l}{1+q})_i}{\Gamma(i+1)\Gamma(m-i+1)}
\frac{1}{\Gamma(1-m) \prod_{l=0}^q (\frac{l+1-m}{1+q})_i}
.
\end{equation}
In \eqref{eq.2sum} the final binomial term may be zero for small $i$ if
$k+(q+1)i< k+m$, which causes spurious singularities in the previous
equation at small $i$. It's advisable to omit these by reversing
the sign of the $i$-sum to start at $m$ instead of 0, $i\to m-i$:
\begin{equation}
a_{k,q}(m)
=
\frac{\Gamma(m+1)}{\Gamma(1+k+m)}
\sum_{i=0}^m 
(-)^m
\frac{\Gamma(1+k) }{\Gamma(i+1)\Gamma(m-i+1)\Gamma(1-m)}
\prod_{l=0}^q 
\frac{
(\frac{k+1+l}{1+q})_{m-i}
}
{
(\frac{l+1-m}{1+q})_{m-i}
}
.
\end{equation}
The substitution \cite[(I.9)]{SlaterHyp}\cite{RoyAMM94}
\begin{equation}
(a)_{N-n}=\frac{(-)^n(a)_N}{(1-a-N)_n}
\end{equation}
yields
\begin{equation}
a_{k,q}(m)
=
\frac{\Gamma(m+1)}{\Gamma(1+k+m)}
\sum_{i=0}^m 
(-)^m
\frac{\Gamma(1+k) }{\Gamma(i+1)\Gamma(m-i+1)\Gamma(1-m)}
\prod_{l=0}^q 
\frac{
(\frac{k+1+l}{1+q})_m
(1-\frac{l+1-m}{1+q}-m)_i
}
{
(\frac{l+1-m}{1+q})_m
(1-\frac{k+1+l}{1+q}-m)_i
}
.
\end{equation}
Reverse application of the $\Gamma$-product formula
\cite{RagabGMJ6}\cite[6.1.20]{AS}
\begin{equation}
\Gamma(z)\Gamma(z+\frac{1}{n})\cdots \Gamma(z+\frac{n-1}{n}) = (2\pi)^{(n-1)/2}n^{1/2-nz}\Gamma(nz); \quad nz\neq 0, -1,-2,\ldots.
\end{equation}
and applying \eqref{eq.negarg} at $n\to m$, $k\to i$ to replace $\Gamma(m+1-i)$
finally renders this as
\begin{multline}
a_{k,q}(m)
=
\binom{m(q+1)+k}{k+m}
\\ \times
{}_{q+2}F_{q+1}
\left(\begin{array}{c}
-m,
1-\frac{0+1-m}{q+1}-m,
1-\frac{1+1-m}{q+1}-m,
\cdots,
1-\frac{q+1-m}{q+1}-m
\\
1-\frac{k+1+0}{q+1}-m,
1-\frac{k+1+1}{q+1}-m,
\cdots,
1-\frac{k+1+q}{q+1}-m
\end{array}\mid 1\right)
.
\end{multline}

\section{Keeping Summation Orders} 
With the same technique as in the previous section --- transforming
the $\Gamma$-products until the index summation variable
 $i$ is the lone index in the Pochhammer symbols --- the
sum \eqref{eq.defb} can be translated to
a terminating Generalized Hypergeometric Series:
\begin{multline}
b_{k,q}(j)=\Gamma(j+1)\sum_{i=0}^j (-)^i\frac{\Gamma(j+k+qi+1)}{\Gamma(i+1)\Gamma(j-i+1)\Gamma(1+j+k)\Gamma(1+qi))}
\\
=
{}_{q+1}F_q
\left(\begin{array}{c}
-j , 
\frac{k+j+1}{q},
\frac{k+j+2}{q},
\cdots
\frac{k+j+q}{q}\\
\frac{1}{q},
\frac{2}{q},
\cdots
\frac{q}{q}
\end{array}\mid 1\right)
.
\end{multline}
Note that this is not equivalent to the last binomial term
in \eqref{eq.2sum}; it does not depend on $m$.
This means that
\begin{equation}
a_{k,q}(m) = \sum_{j=0}^m (-)^j \binom{m}{j} b_{k,q}(j)
\label{eq.abinb}
\end{equation}
represents $a$ as a binomial transform  of $b$, whereas \eqref{eq.2sum}
is \emph{not} a binomial transform pair.

A special case is the Gaussian Hypergeometric Function \cite[15.4.24]{DLMF}
\begin{equation}
b_{k,1}(j) = \binom{-k-1}{j}
\end{equation}
with the obvious recurrence
\begin{equation}
(j+1)b_{k,1}(j+1)
+
(k+j+1)b_{k,1}(j)
=0.
\end{equation}

The generating function is defined as
\begin{equation}
B_{k,q}(z) \equiv \sum_{j\ge 0} b_{k,q}(j)z^j,
\end{equation}

with examples in Table \ref{tab.bkq}.
\begin{table}
\begin{tabular}{rr|l|l}
$k$ & $q$ & $b_{k,q}(j)$ & $B_{k,q}(z)$ \\
\hline
0 & 0 & 1 0 0 0 0 0 0 0 0 0 0 0 0 0 0 0& 1 \\
0 & 1 & 1 -1 1 -1 1 -1 1 -1 1 -1 1 -1 1 -1 1 -1 & $\frac{1}{1+z}$ \\
0 & 2 & 1 -2 4 -8 16 -32 64 -128 256 -512 1024 -2048  & $\frac{1}{1+2z}$ \\
0 & 3 & 1 -3 9 -27 81 -243 729 -2187 6561 -19683 59049 -177147 & $\frac{1}{1+3z}$ \\
1 & 0 & 1 0 0 0 0 0 0 0 0 0 0 0 0 0 0 0 & $1$ \\
1 & $1/2$ & 1 -7/8 5/8 -13/32 1/4 -19/128 11/128 & $\frac{8+z}{2(z+2)^2}$ \\
1 & 1 & 1 -2 3 -4 5 -6 7 -8 9 -10 11 -12 13 -14 15 -16 & $\frac{1}{(1+z)^2}$ \\
1 & $3/2$ &  1 -27/8 63/8 -513/32 243/8 -7047/128 12393/128 -85293/512 & $\frac{8-3z}{2(2+3z)^2}$\\
1 & 2 & 1 -5 16 -44 112 -272 640 -1472 3328 -7424 16384 -35840& $\frac{1-z}{(1+2z)^2}$\\
1 & 3 & 1 -9 45 -189 729 -2673 9477 -32805 111537 -373977 1240029 -4074381 & $\frac{1-3z}{(1+3z)^2}$\\
1 & 4 & 1 -14 96 -544 2816 -13824 65536 -303104 1376256 -6160384 27262976 & $\frac{1-6z}{(1+4z)^2}$ \\
1 & 5 & 1 -20 175 -1250 8125 -50000 296875 -1718750 9765625 -54687500 & $\frac{1-10z}{(1+5z)^2}$ \\
2 & 0 & 1 0 0 0 0 0 0 0 0 0 0 0 0 0 0 0 & $1$\\
2 & 1 & 1 -3 6 -10 15 -21 28 -36 45 -55 66 -78 91 -105 120 -136& $\frac{1}{(1+z)^3}$\\
2 & 2 & 1 -9 41 -146 456 -1312 3568 -9312 23552 -58112 140544 & $\frac{1-3z-z^2}{(1+2z)^3}$ \\
2 & 3 & 1 -19 141 -783 3753 -16443 67797 -267543 1021329 -3798819 & $\frac{1-10z-3z^2}{(1+3z)^3}$ \\
2 & 4 & 1 -34 356 -2704 17536 -103424 572416 -3026944 15466496 -76939264 & $\frac{1-22z-4z^2}{(1+4z)^3}$ \\
2 & 5 & 1 -55 750 -7250 59375 -440625 3062500 -20312500 130078125 & $\frac{1-40z}{(1+5z)^3}$ \\
3 & 0 & 1 0 0 0 0 0 0 0 0 0 0 0 0 0 0 0& $1$ \\
3 & 1 & 1 -4 10 -20 35 -56 84 -120 165 -220 286 -364 455 -560 680 -816 &$\frac{1}{(1+z)^4}$ \\
3 & 2 & 1 -14 85 -377 1408 -4712 14608 -42800 120064 -325376 857344  &$\frac{1-6z-3z^2-z^3}{(1+2z)^4}$ \\
3 & 3 & 1 -34 351 -2484 14445 -74358 352107 -1568808 6672537 -27359370  &$\frac{1-22z-3z^2}{(1+3z)^4}$ \\
3 & 4 & 1 -69 1036 -10184 80896 -564224 3603456 -21592064 123273216  &$\frac{(1-z)(1-52z-24z^2)}{(1+4z)^4}$ \\
\end{tabular}
\caption{Examples of $b_{k,q}(j)$, $j\ge 0$,
and 
generating functions for small indices.}
\label{tab.bkq}
\end{table}

For $k=0$ we have in particular
\begin{equation}
b_{0,q}(j) = (-q)^j
\label{eq.b0q}
\end{equation}
and consequently the geometric series
\begin{equation}
B_{0,q}(z)=\sum_{j\ge 0} (-qz)^j = \frac{1}{1+qz}.
\label{eq.B0qz}
\end{equation}
\begin{proof}
The ratio of consecutive terms of $j$ and $j+1$ for $k=0$ in \eqref{eq.defb} is by standard recurrences of binomials
\begin{equation}
\frac{ \binom{j+1}{i}\binom{j+1+qi}{j+1}}
{\binom{j}{i}\binom{j+qi}{j}}
= \frac{j+1+iq}{j+1-i}
= -q+\frac{(j+1)(q+1)}{j+1-i}
\end{equation}
and therefore
\begin{equation}
(-)^i\binom{j+1}{i}\binom{j+1+iq}{j+1}
=
-q(-)^i \binom{j}{i}\binom{j+iq}{j}
+(q+1) (-)^i \binom{j}{i}\binom{j+iq}{j}\frac{j+1}{j+1-i}
\end{equation}
\begin{equation}
=
-q (-)^i \binom{j}{i}\binom{j+iq}{j}
+(q+1) (-)^i \binom{j+1}{i}\binom{j+iq}{j}.
\end{equation}
Summation over $i$ (where the first term on the right hand side is zero at the upper limit $i=j+1$) yields
\begin{equation}
\sum_{i=0}^{j+1} (-)^i \binom{j+1}{i}\binom{j+1+iq}{j+1}
=
-q\sum_{i=0}^{j(+1)} (-)^i \binom{j}{i}\binom{j+iq}{j}
+(q+1)\sum_{i=0}^{j+1} (-)^i \binom{j+1}{i}\binom{j+iq}{j}
.
\label{eq.b0qtmp}
\end{equation}
We show that the last term on the right hand side is zero by transforming the sum
to a (Saalsch\"utzian) terminating hypergeometric series in the same technique as above:
\begin{equation}
\sum_{i=0}^{j+1} (-)^i \binom{j+1}{i}\binom{j+iq}{j} 
=
{}_qF_{q+1}\left(\begin{array}{c}
-j-1, \frac{j+1}{q}, \frac{j+2}{q}, \cdots ,\frac{j+q}{q} \\
\frac{1}{q}, \frac{2}{q}, \cdots ,\frac{q}{q} 
\end{array}\mid 1\right)
\end{equation}
and rephrase this as a finite integral \cite{CoffeyJCAM233,DriverETNA25}
\begin{equation}
=
\frac{\Gamma(1)}{\Gamma(j+1)\Gamma(-j)}
\int_0^1 t^j (1-t)^{-j-1}(1-t^q)^{j+1}dt
=
\frac{\Gamma(1)}{\Gamma(j+1)\Gamma(-j)}
\int_0^1 t^j (\frac{1-t^q}{1-t})^{j+1}dt
=0.
\end{equation}
The last equation follows because  \cite[24.1.2]{AS}
\begin{multline}
\int_0^1 t^j (\frac{1-t^q}{1-t})^{j+1}dt
=\int_0^1 t^j (1+t+t^2+\cdots t^{q-1})^{j+1}dt
\\
=\sum_{j_0+j_1+j_2+\cdots j_{q-1}=j+1}
\binom{j+1}{j_0,j_1,j_2,\ldots,j_{q-1}}
\frac{1}{1+j_1+2j_2+3j_3+\cdots (q-1)j_{q-1}}
\end{multline}
integrates to a positive finite rational number,
and because the pole of $\Gamma(-j)$ in the denominator grows to infinity.
Therefore \eqref{eq.b0qtmp} is effectively
\begin{equation}
b_{0,q}(j+1)
=
-qb_{0,q}(j)
\end{equation}
and solving this recurrence starting at the obvious $b_{0,q}(0)=1$ produces \eqref{eq.b0q}.
\end{proof}
The generating functions for $k>0$ are recursively obtained via
\begin{equation}
B_{k,q}(z)
=
B_{k-1,q}(z)
+
\frac{1}{(1+qz)^{k+1}}
\sum_{s=0}^{k} 
z^{s}
\sum_{j=0}^{s}
\binom{k+1}{s-j}q^{s-j}
\sum_{i=0}^j (-)^i\binom{j}{i}\binom{j+k-1+qi}{qi-1}
\label{eq.Bkqzrec}
\end{equation}
starting at \eqref{eq.B0qz}---the derivation is skipped in this manuscript. 
It shows that the generating functions are polynomials
in $z$ divided through $(1+qz)^{k+1}$. The $b_{k,q}(j)$ are C-finite
functions.

\section{GF of the Main Sequence}
If
\begin{equation}
b_n=\sum_{k=0}^n\binom{n}{k}a_k \leftrightarrow a_n=\sum_{k=0}^n(-1)^{n-k}\binom{n}{k}b_k,
\label{eq.bindef}
\end{equation}
is a binomial transform pair, then the action on the GF's is to substitute
the arguments by rational functions
\cite{SprugnoliDM142,BernsteinLAA226}
\begin{equation}
B(x) = \frac{1}{1-x}A(\frac{x}{1-x})
\leftrightarrow
A(x) = \frac{1}{1+x}B(\frac{x}{1+x})
\end{equation}
This mapping applies to \eqref{eq.abinb}:
\begin{defn}
 (ordinary generating function)
\begin{equation}
A_{k,q}(z)\equiv \sum_{m\ge 0} a_{k,q}(m) z^m.
\label{eq.Akqdef}
\end{equation}
\end{defn}
The C-finite GF's of \eqref{eq.Bkqzrec} transform into C-finite GF's:
\begin{equation}
A_{k,q}(z)=\frac{1}{1-z} B_{k,q}\left(\frac{-z}{1-z}\right).
\label{eq.AfromB}
\end{equation}
(The switch $z\to -z$ is an adaptation to the sign choices in \eqref{eq.bindef}).
Table \ref{tab.Aofz} illustrates this for small $k$ and $q$.

\begin{table}
\begin{tabular}{rr|rrrrrrrrr|l}
$k$ & $q$ & $A_{k,q}$ \\
\hline
0 & 0 & $-1/(-1+z)$\\
0 & 1 & $-1/(-1+2 z)$\\
0 & 2 & $-1/(-1+3 z)$\\
0 & 3 & $-1/(-1+4 z)$\\
0 & 4 & $-1/(-1+5 z)$\\
0 & 5 & $-1/(-1+6 z)$\\
1 & 0 & $-1/(-1+z)$\\
1 & 1 & $-(-1+z)/(-1+2 z)^2$\\
1 & 2 & $1/(-1+3 z)^2$ & A027471\\
1 & 3 & $(1+2 z)/(-1+4 z)^2$\\
1 & 4 & $(1+5 z)/(-1+5 z)^2$\\
1 & 5 & $(1+9 z)/(-1+6 z)^2$\\
2 & 0 & $-1/(-1+z)$\\
2 & 1 & $-(-1+z)^2/(-1+2 z)^3$\\
2 & 2 & $(-1-z+3 z^2)/(-1+3 z)^3$\\
2 & 3 & $(-1-8 z+12 z^2)/(-1+4 z)^3$ & A361609\\
2 & 4 & $(-1-20 z+25 z^2)/(-1+5 z)^3$\\
2 & 5 & $(-1+z) (1+39 z)/(-1+6 z)^3$\\
3 & 0 & $-1/(-1+z)$\\
3 & 1 & $-(-1+z)^3/(-1+2 z)^4$\\
3 & 2 & $(1+3 z-12 z^2+9 z^3)/(-1+3 z)^4$\\
3 & 3 & $(-1+z) (-1-20 z+24 z^2)/(-1+4 z)^4$\\
3 & 4 & $-(75 z^2-50 z-1)/(-1+5 z)^4$ & A361610\\
3 & 5 & $-(-1-102 z+57 z^2+171 z^3)/(-1+6 z)^4$\\
4 & 0 & $-1/(-1+z)$\\
4 & 1 & $-(-1+z)^4/(-1+2 z)^5$\\
4 & 2 & $(-1-6 z+29 z^2-39 z^3+18 z^4)/(-1+3 z)^5$\\
4 & 3 & $-(1+36 z-92 z^2+48 z^3+16 z^4)/(-1+4 z)^5$\\
4 & 4 & $-(1+101 z-65 z^2-425 z^3+500 z^4)/(-1+5 z)^5$\\
4 & 5 & $-(1+222 z+388 z^2-2496 z^3+2385 z^4)/(-1+6 z)^5$\\
5 & 0 & $-1/(-1+z)$\\
5 & 1 & $-(-1+z)^5/(-1+2 z)^6$\\
5 & 2 & $(1+10 z-55 z^2+99 z^3-81 z^4+27 z^5)/(-1+3 z)^6$\\
5 & 3 & $-(-1-60 z+132 z^2+100 z^3-432 z^4+288 z^5)/(-1+4 z)^6$\\
5 & 4 & $-(-1-180 z-270 z^2+2800 z^3-4625 z^4+2500 z^5)/(-1+5 z)^6$\\
5 & 5 & $-(-1+z) (1+427 z+3123 z^2-10206 z^3+7155 z^4)/(-1+6 z)^6$\\
5 & 6 & $(36015 z^4-40474 z^3+10731 z^2+882 z+1)/(-1+7 z)^6$ & A361608\\
\end{tabular}
\caption{Examples of generating functions $A_{k,q}(z)$ for small indices and a few associated sequence numbers of the
Online Encyclopedia of Integer Sequences \cite{sloane}.}
\label{tab.Aofz}
\end{table}

\section{Summary}
The ordinary generating function \eqref{eq.Akqdef}
is obtained starting from \eqref{eq.B0qz}, creating generating
functions for larger $k$ recursively via \eqref{eq.Bkqzrec}, and finally
reflecting the arguments via \eqref{eq.AfromB}.

\appendix
\section{GF of the Aerated Binomial Term}
The core ingredient of \eqref{eq.defb} is
\begin{defn} (aerated binomial)
\begin{equation}
c_{J,q}(i)
\equiv \binom{J+qi}{J},
\end{equation}
\end{defn}
rewritten as products of Pochhammer symbols with \eqref{eq.pmulti}
\begin{equation}
c_{J,q}(i)=\frac{\Gamma(J+qi+1)}{\Gamma(1+J)\Gamma(qi+1)} 
= \frac{(J+1)_{qi}}{(1)_{qi}}
= \frac{((J+1)/q)_i ((J+2)/q)_i\cdots ((J+q)/q)_i }{ (1/q)_i (2/q)_i\cdots (q/q)_i}
.
\end{equation}
Its generating function is
\begin{equation}
C_{J,q}(z) \equiv \sum_{i\ge 0} c_{J,q}(i)z^i
=
{}
_{q+1}F_q(\begin{array}{c}
1, \frac{J+1}{q}, \frac{J+2}{q},\cdots ,\frac{J+q}{q}\\
\frac{1}{q}, \frac{2}{q},\cdots ,\frac{q}{q}
\end{array}\mid z)
=
{}
_{q}F_{q-1}(\begin{array}{c}
\frac{J+1}{q}, \frac{J+2}{q},\cdots ,\frac{J+q}{q}\\
\frac{1}{q}, \frac{2}{q},\cdots ,\frac{q-1}{q}
\end{array}\mid z)
\end{equation}
\begin{exa}
\begin{equation}
c_{J,0}(i) = c_{0,q}(i) = 1;\quad C_{J,0}(z) = C_{0,q}(z) = \frac{1}{1-z}.
\end{equation}
\end{exa}
\begin{exa}
\begin{equation}
c_{J,1}(i) = \binom{J+i}{J};\quad C_{J,1}(z) 
= 
{}
_{2}F_1(\begin{array}{c}
1, J+1\\
1
\end{array}\mid z)
= 
{}
_1F_0(\begin{array}{c}
J+1\\
\end{array}\mid z)
=\frac{1}{(1-z)^{J+1}}
\end{equation}
\end{exa}
\begin{exa}
An application of the Eulerian polynomials is \cite{NealCMJ25}
\begin{equation}
c_{1,q}(i) = 1+qi ;\quad C_{q,z}(z) = \sum_{i\ge 0}z^i + q\sum_{i\ge 0 }iz^i = 
\frac{1+(q-1)z}{(1-z)^2}.
\end{equation}
\end{exa}

\begin{exa}
\begin{equation}
c_{J,2}(i) = \binom{J+2i}{J};\quad C_{J,2}(z) 
= 
{}
_2F_1(\begin{array}{c}
\frac{J+1}{2},\frac{J+2}{2}\\
\frac{1}{2}
\end{array}\mid z)
\label{eq.CJ2z}
\end{equation}
The basic values at $J=0$ and $1$ are
\begin{eqnarray}
c_{0,2}(i)=1&,& \quad C_{0,2}(z)=\frac{1}{1-z};\\
c_{1,2}(i)=1+2i&,& \quad C_{1,2}(z)=\frac{1+z}{(1-z)^2}.
\end{eqnarray}
A contiguous relation of the Gaussian Hypergeometric Function \cite[15.5.13]{DLMF}
is
\begin{equation}
(c-a-b-1)F(a,b+1;c;z)+a(1-z)F(a+1,b+1;c;z)-(c-b-1)F(a,b;c;z)=0.
\end{equation}
Taking $b=(J+1)/2$, $a=J/2+1$, $c=1/2$ and translation
with \eqref{eq.CJ2z} yields the recurrence
\begin{equation}
 (1-z)C_{J,2}(z) = 2C_{J-1,2}(z)- C_{J-2,2}(z).
\end{equation}
which can be summarized as rational generating functions
\cite[A034839]{sloane}
\begin{equation}
C_{J,2}(z) = \frac{\sum_{l=0}^{\lfloor L/2\rfloor}\binom{J+1}{2l}z^l}{(1-z)^{J+1}}
\end{equation}
\end{exa}

\begin{table}
\begin{tabular}{rr|rrrrrrrrr|l}
$J$ & $q$ & 1 & 2& 3& 4& 5& 6 & 7 & 8& 9& $C_{J,q}$ \\
\hline
1 & 3 &1&4&7&10&13&16&19&22&25 & $\frac{1+2x}{(1-x)^2}$\\
1 & 4 &1&5&9&13&17&21&25&29&33 & $\frac{1+3x}{(1-x)^2}$\\
1 & 5 &1&6&11&16&21&26&31&36&41 & $\frac{1+4x}{(1-x)^2}$\\
2 & 3 &1&10&28&55&91&136&190&253&325 & $\frac{1+7x+x^2}{(1-x)^3}$\\
2 & 4 &1&15&45&91&153&231&325&435&561 & $\frac{1+12x+3x^2}{(1-x)^3}$\\
2 & 5 &1&21&66&136&231&351&496&666&861 & $\frac{1+18x+6x^2}{(1-x)^3}$\\
3 & 3 &1&20&84&220&455&816&1330&2024&2925 & $\frac{1+16x+10x^2}{(1-x)^4}$\\
3 & 4 &1&35&165&455&969&1771&2925&4495&6545 & $\frac{(1+x)(1+30x+x^2)}{(1-x)^4}$\\
3 & 5 &1&56&286&816&1771&3276&5456&8436&12341 & $\frac{1+52x+68x^2+4x^3}{(1-x)^4}$\\
4 & 3 &1&35&210&715&1820&3876&7315&12650&20475 & $\frac{1+30x+45x^2+5x^3}{(1-x)^5}$\\
4 & 4 &1&70&495&1820&4845&10626&20475&35960&58905 & $\frac{1+65x+155x^2+35x^3}{(1-x)^5}$\\
4 & 5 &1&126&1001&3876&10626&23751&46376&82251&135751 & $\frac{1+121x+381x^2+121x^3+x^4}{(1-x)^5}$\\
5 & 3 &1&56&462&2002&6188&15504&33649&65780&118755 & $\frac{1+50x+141x^2+50x^3+x^4}{(1-x)^6}$\\
5 & 4 &1&126&1287&6188&20349&53130&118755&237336&435897 & $\frac{1+120x+546x^2+336x^3+21x^4}{(1-x)^6}$\\
5 & 5 &1&252&3003&15504&53130&142506&324632&658008&1221759 & $\frac{1+246x+1506x^2+1246x^3+126x^4}{(1-x)^6}$\\
\end{tabular}
\caption{Examples of $c_{J,q}(i)$ for small indices and their generating functions $C_{J,q}(x)$, see \eqref{eq.Cofz}.}
\end{table}

Perfect powers have representations
\cite{SinghJIS19}
\begin{equation}
k^n=
\sum_{j=0}^n \begin{Bmatrix}n\\j\end{Bmatrix} \binom{k}{j}j!,
\end{equation}
where $\begin{Bmatrix} .\\.\end{Bmatrix}$ are the Stirling Numbers of the Second Kind.
Taking the generating function with respect to $k$, transforming the emergent
hypergeometic series with one of the linear transformations yields
rational generating functions
\cite{PizaMM21} 
\cite[(5.2)]{GouldFQ16}
\cite{KrishnapCMJ26}
\cite[(24)]{Boyadzhievarxiv2011}
\cite{BoyadAM42}
\begin{equation}
\sum_{k=0}^\infty k^nx^k = \frac{1}{1-x}\omega_n(\frac{x}{1-x})
=\frac{\sum_{k=0}^n \left\langle \begin{matrix}n\\k \end{matrix}\right\rangle x^k}{(1-x)^{n+1}}
.
\end{equation}
The coefficients $\left\langle \begin{matrix}.\\. \end{matrix}\right\rangle$ of
the polynomial in the numerator are the Eulerian Numbers \cite[A123125]{sloane}.
\begin{defn} (geometric polynomials)
\begin{equation}
\omega_n(x)=\sum_{k=0}^n \begin{Bmatrix} n \\ k\end{Bmatrix} k!x^k
.
\label{eq.omegadef}
\end{equation}
\end{defn}
\begin{exa}
\begin{eqnarray}
\omega_1(x)&=&x;\\
\omega_2(x)&=&x+x^2;\\
\omega_3(x)&=&x+6x^2+6x^3;\\
\omega_4(x)&=&x+14x^2+36x^3+24x^4.
\end{eqnarray}
\end{exa}
The coefficients $[x^k]\omega_n(x)$ are tabulated 
in the OEIS \cite[A019538]{sloane}. 
\begin{rem}
The generating function of the $\omega$-polynomials is
\cite[(3.4)]{BoyadIJMMS2005}
\begin{equation}
x^n=\frac{1}{n!}\sum_{k=0}^n \begin{bmatrix} n\\ k\end{bmatrix}\omega_k(x),
\end{equation}
where $\begin{bmatrix} .\\. \end{bmatrix}$ are the (signed) Stirling numbers of the First Kind \cite[24.1.3]{AS}\cite[26.8]{DLMF}.
This is the generic Stirling transform of \eqref{eq.omegadef} using the orthogonality of the Stirling numbers
of both kinds.
\end{rem}

With the standard generating functions for the Stirling Numbers of the First Kind \cite[Table 264]{GrahamKnuthPatashnik}, our
atomic binomial becomes
\begin{multline}
c_{J,q}(i)=\binom{J+iq}{J}
=\frac{(J+iq)(J+iq-1)\cdots (iq+1)}{1\times 2\times \cdots \times J}
=\frac{1}{J!}\sum_{l=0}^J \begin{bmatrix} J\\ l\end{bmatrix} (J+iq)^l
\\
=\frac{1}{J!}\sum_{l=0}^J \begin{bmatrix} J\\ l\end{bmatrix}\sum_{t=0}^l \binom{l}{t}(iq)^tJ^{l-t}
\\
=\frac{1}{J!}\sum_{t=0}^J \sum_{l=t}^J \begin{bmatrix} J\\ l\end{bmatrix}\binom{l}{t}q^t i^t J^{l-t}
.
\end{multline}
Its generating functions are
\begin{multline}
C_{J,q}(z)
=
\frac{1}{J!}\sum_{i\ge 0}\sum_{t=0}^J \sum_{l=t}^J \begin{bmatrix} J\\ l\end{bmatrix} \binom{l}{t}
q^tJ^{l-t}i^tz^i
=
\frac{1}{J!}\sum_{t=0}^J \sum_{l=t}^J \begin{bmatrix} J\\ l\end{bmatrix} \binom{l}{t}
q^tJ^{l-t}
\frac{1}{1-z}\omega_t(\frac{z}{1-z})
\\
=
\frac{1}{J!}\sum_{t=0}^J 
\frac{1}{1-z}\omega_t(\frac{z}{1-z})
q^t
\sum_{l=t}^J \begin{bmatrix} J\\ l\end{bmatrix} \binom{l}{t}J^{l-t}
\\
=
\frac{1}{J!}\sum_{t=0}^J 
\frac{1}{1-z}\omega_t(\frac{z}{1-z})
q^t
(-)^{J+t}
\begin{bmatrix}
J+1 \\ t+1
\end{bmatrix}.
\label{eq.Cofz}
\end{multline}

\begin{rem}
The transition to the last equation is a partial binomial transform
\begin{equation}
\sum_{l=t}^J \begin{bmatrix} J\\l \end{bmatrix} \binom{l}{t}J^l= (-)^{J+t}J^t \begin{bmatrix} J+1\\ t+1\end{bmatrix}
.
\end{equation}

This is an associate of \cite[(6.16)]{GrahamKnuthPatashnik}
\begin{equation}
\sum_{k=0}^n (-)^{n+k} \begin{bmatrix} n\\k \end{bmatrix}\binom{k}{m}=(-)^{n+m} \begin{bmatrix}n+1\\ m+1\end{bmatrix},
\end{equation}
once the rising $k$-binomial transform of a sequence $a_i$ is defined as
\cite{SpiveyJIS9}
\begin{equation}
r_n\equiv \sum_{i=0}^n \binom{n}{i}k^ia_i,
\end{equation}
Spivey's Theorem \cite[thm.3.2]{SpiveyJIS9} is applied to unfold this
as $k$ successive binomial transforms,
and \cite[Table 2.1]{Riordan}
\begin{equation}
b_n=\sum_{k\ge n}\binom{k}{n}a_k \leftrightarrow a_n=\sum_{k\ge n}(-)^{k+n}\binom{k}{n}b_k
\end{equation}
wraps and unwraps these successive binomial transforms.
\end{rem}

\bibliographystyle{amsplain}
\bibliography{all}

\end{document}